\documentclass [12pt]{article}
\usepackage{amsmath,amssymb,amscd}

\begin{document}
\def\l{\lambda}
\def\m{\mu}
\def\a{\alpha}
\def\b{\beta}
\def\g{\gamma}
\def\G{\Gamma}
\def\d{\delta}
\def\e{\epsilon}
\def\o{\omega}
\def\O{\Omega}
\def\v{\varphi}
\def\t{\theta}
\def\r{\rho}
\def\bs{$\blacksquare$}
\def\bp{\begin{proposition}}
\def\ep{\end{proposition}}
\def\bt{\begin{theo}}
\def\et{\end{theo}}
\def\be{\begin{equation}}
\def\ee{\end{equation}}
\def\bl{\begin{lemma}}
\def\el{\end{lemma}}
\def\bc{\begin{corollary}}
\def\ec{\end{corollary}}
\def\pr{\noindent{\bf Proof: }}
\def\note{\noindent{\bf Note. }}
\def\bd{\begin{definition}}
\def\ed{\end{definition}}
\def\C{{\mathbb C}}
\def\P{{\mathbb P}}
\def\Z{{\mathbb Z}}
\def\d{{\rm d}}
\def\deg{{\rm deg\,}}
\def\deg{{\rm deg\,}}
\def\arg{{\rm arg\,}}
\def\min{{\rm min\,}}
\def\max{{\rm max\,}}
 \def\d{{\rm d}}
\def\deg{{\rm deg\,}}
\def\deg{{\rm deg\,}}
\def\arg{{\rm arg\,}}
\def\min{{\rm min\,}}
\def\max{{\rm max\,}}

\newtheorem{theo}{Theorem}[section]
\newtheorem{lemma}{Lemma}[section]
\newtheorem{definition}{Definition}[section]
\newtheorem{corollary}{Corollary}[section]
\newtheorem{proposition}{Proposition}[section]

\begin{titlepage}
\begin{center}

\topskip 5mm

{\LARGE{\bf {Signal Acquisition from Measurements via Non-Linear
Models}}}

\vskip 6mm

{\large {\bf N. Sarig $^{*}$,}} \hspace {2 mm} {\large {\bf Y.
Yomdin$^{*}$}}

\vspace{3 mm}
\end{center}

{$^{*}$ Department of Mathematics, The Weizmann Institute of
Science, Rehovot 76100, Israel}

\vspace{3 mm}

{e-mail: niv.sarig@weizmann.ac.il, \ yosef.yomdin@weizmann.ac.il}

\vspace{3 mm}
\begin{center}

{ \bf Abstract}
\end{center}
{\small{We consider the problem of reconstruction of a non-linear
finite-parametric model $M=M_p(x),$ with $p=(p_1,\dots,p_r)$ a set
of parameters, from a set of measurements $\mu_j(M)$. In this
paper $\mu_j(M)$ are always the moments $m_j(M)=\int x^jM_p(x)dx$.
This problem is a central one in Signal Processing, Statistics,
and in many other applications.

We concentrate on a direct (and somewhat ``naive") approach to the
above problem: we simply substitute the model function $M_p(x)$
into the measurements $\mu_j$ and compute explicitly the resulting
``symbolic" expressions of $\mu_j(M_p)$ in terms of the parameters
$p$. Equating these ``symbolic" expressions to the actual
measurement results, we produce a system of nonlinear equations on
the parameters $p$, which we consequently try to solve.

The aim of this paper is to review some recent results (mostly of
\cite{Vet5,Mil1,Mil2,Put1,Vet4,Vet3,Mil3,Mil4,Vet2}) in this
direction, stressing the algebraic structure of the arising
systems and mathematical tools required for their solutions.

In particular, we discuss the relation of the reconstruction
problem above with the recent results of
\cite{bfy,bry,chr,pak1,pak2,pak3,pry,ry} on the vanishing problem
of generalized polynomial moments and on the Cauchy-type integrals
of algebraic functions.

The accompanying paper \cite{Kis1} (this volume) provides a
solution method for a wide class of reconstruction problems as
above, based on the study of linear differential equations with
rational coefficient, which are satisfied by the moment generating
function of the problem.

\vspace{2 mm}
\begin{center}
------------------------------------------------
\vspace{2 mm}
\end{center}
This research was supported by the ISF, Grant No. 264/05 and by
the Minerva foundation.

\smallskip

Key words: Signal acquisition, Non-linear models, Moments
inversion.

\smallskip

AMS Subject Classification: 62J02, 14P10, 42C99}}

\end{titlepage}
\newpage

%SECTION 1
%%%%%%%%%%%%%%%%%%%%%%%%%%%%%%%%%%%%%%%%%%%%%%%%%%%%%%%%%%%%%%%%%%%%%%%%%%%%%%%

\section{Introduction}
\setcounter{equation}{0}

In this paper we consider the following problem: let a
finite-parametric family of functions $M=M_p(x), \ x\in {\mathbb
R}^m$ be given, with $p=(p_1,\dots,p_r)$ a set of parameters. We
call $M_p(x)$ a model, and usually we assume that it depends on
some of its parameters in a non-linear way (this is always the
case with the ``geometric" parameters representing the shape and
the position of the model).

The problem is:

\smallskip

\noindent{\it How to reconstruct in a robust and efficient way the
parameters $p$ from a set of ``measurements"
$\mu_1(M),\dots,\mu_l(M)$?}

\smallskip

In this paper $\mu_j$ will be the moments $m_j(M)=\int
x^jM_p(x)dx$. This assumption is not too restrictive - see, for
example, \cite{Mil1,Vet5}.

The above problem is certainly among the central ones in Signal
Processing (Non-linear matching), Statistics (Non-linear
Regression), and in many other applications. See
\cite{Vet5,Mil1,Mil2,Vet4,Vet3,Mil3,Mil4,Vet2} and references
there.

We concentrate in the present paper on a direct (and somewhat
``naive") approach to the above problem: we simply substitute the
model function $M_p(x)$ into the measurements $\mu_j$ and compute
explicitly the resulting ``symbolic" expressions of $\mu_j(M_p)$
in terms of the parameters $p$. Equating these ``symbolic"
expressions to the actual measurement results, we produce a system
of nonlinear equations on the parameters $p$ which we try to
solve.

\smallskip

Certainly, the polynomial moments do not present the best choice
of measurements for practical applications since the monomials
$x^j$ are far away from being orthogonal (see, for example,
\cite{Tal}). However, the main features of the arising non-linear
systems remain the same for a much wider class of measurements,
while their structure is much more transparent for moments.

\smallskip

The aim of this paper is to review some recent results (mostly of
\cite{Vet5,Mil1,Mil2,Put1,Vet4,Vet3,Mil3,Mil4,Vet2}) in this
direction, stressing the algebraic structure of the arising
systems and mathematical tools required for their solutions. In
particular, we stress the role of the moment generating function.

We start with some initial examples of the models $M_p(x)$ in one
dimension: these are polynomials and rational functions. Then we
consider linear combinations of $\delta$-functions. The system
which appears in this example is typical in many application. We
discuss one of the classical solutions methods, following
\cite{Pro,Kis2,Mil1,Mil2,Vet4,Vet2}.

Next we deal with piecewise-solution of linear differential
equations, providing some pre-requisites for the reconstruction
method described in \cite{Kis1} (this volume). Then we consider
piecewise-algebraic functions of one variable. We prove
injectivity of the finite moment transform on such functions, and
discuss the relation of the reconstruction problem for such
functions with the recent results of
\cite{bfy,bry,bry1,chr,pak1,pak2,pak3} on the vanishing problem of
generalized polynomial moments.

In two dimensions we shortly present results of
\cite{Mil1,Mil2,Mil3} concerning reconstruction of polygons from
their complex moments, as well as results of \cite{Put1} on
reconstruction of ``quadrature domains". Finally we consider the
problem of reconstruction of $\delta$-functions along algebraic
curves, relating it to the vanishing problem of double moments
(\cite{Wer,Hen1,Hen2,Hen3,bry,pry}).

We almost do not touch the classical Moment Theory, refereeing the
reader to \cite{Nik.Sor} and especially to
\cite{Put1,Put2,Put3,Put4,Put5,Put6,Tal}, where, in particular, a
review of the classical results and methods is given, as applied
to the effective reconstruction problem.

We also don't discuss in this paper the problem of noise
resistance. It is treated in \cite{Mil1,Mil2,Vet4,Vet3,Mil4}.

\subsection{Applicability of the ``direct substitution" method}

The key condition for applicability of our approach is the
assumption that the signals we work with can be faithfully
approximates by a priori known ``simple" geometric models.

A natural question is: to what extent this assumption is
realistic? The answer to this question is twofold:

\medskip

1. In many specific application the form of the signal is indeed
known a priori. Besides the wide circle of applications mentioned
in \cite{Vet5,Mil1,Mil2,Vet4,Vet3,Mil3,Mil4,Vet2} notice that this
is usually the case in visual quality inspection. Similar
situations arise in some medical applications where a non-linear
parametric model of an important pattern has to be matched to the
radiology or ultrasound measurements data.

\medskip

2. A general applicability of our approach in problems involving
image acquisition, analysis and processing depends on a
possibility to represent general images of the real world via
geometric models.

The importance of such a representation in many imaging problems,
from still and video-compression to visual search and pattern
detection is well-recognized. Some initial implementations of
geometric image ``modelization" have been suggested, in
particular, in \cite{Kun.Iko.Koc,Eli.Yom,Bri.Eli.Yom,Eld}. See
\cite{Eld} and references there for a general overview and
analysis of the performance of edges-based methods in images
representation.

However, in general the ``geometric" methods, as for today, suffer
from an inability to achieve a full visual quality for high
resolution photo-realistic images of the real world. {\it In fact,
the mere possibility of a faithful capturing such images with
geometric models presents one of important open problems in Image
Processing, sometimes called ``the vectorization problem"}.

Certainly, this current status of affairs makes problematic
immediate practical applications of general imaging methods based
on geometric model.

Let us express our strong belief that a full visual quality
geometric-model representation of high resolution photo-realistic
images is possible. As achieved, it promises a major advance in
image compression, capturing, and processing, in particular, via
the approach of the present paper.

\medskip

Recently some ``semi-linear" approaches have emerged providing a
reliable reconstruction of ``simple" (and not necessarily regular)
signals from a small number of measurements. In these approaches
(see \cite{Can,Don} and references there) ``simplicity" or
``compressibility" of a function is understood as a possibility
for its accurate sparse representation in a certain (wavelet)
basis.

A somewhat more general approach to the notion of a complexity of
functions has been suggested in \cite{Yom1,Yom2}: here we take as
a complexity measure the rate of semi-algebraic approximation. If
the wavelet base is semi-algebraic, ``compressible" functions have
low semi-algebraic complexity. The same is typically true for
functions allowing for a fast approximation by various types of
non-linear models.

\section{Examples of moment inversion: one variable}
\setcounter{equation}{0}

In this section we consider some natural examples of the models
$M_p(x)$ in one dimension and of their reconstruction from the
moments. These are polynomials, rational functions, linear
combinations of $\delta$-functions, and the class $A_D$ of
piecewise-analytic functions, each piece satisfying a fixed linear
differential operator $D$ with rational coefficients.
(Piecewise-polynomials belong to $A_D$ for $D={{d^n}\over
{dx^n}}$). Then we consider piecewise-algebraic functions.

\medskip

In this paper we use as one of the main tools in solving the
moment inversion problem the moment generating function $I_g(z)$
defined as \be I_g(z)=\sum_{k=0}^{\infty} m_k(g) z^k=\int_0^1
{{g(t)dt}\over {1-zt}}. \ee

\subsection{Vetterli's approach}

In \cite{Vet2,Vet4,Vet5} an important class of signals has been
introduced, possessing a ``finite rate of innovation", i.e. a
finite number of degrees of freedom per unit of time. Usually such
signals are not band-limited, so classical sampling theory does
not enable a perfect reconstruction of signals of this type. In
\cite{Vet2,Vet4,Vet5} it was shown that using an adequate sampling
kernel and a sampling rate greater or equal to the rate of
innovation, it is possible to reconstruct such signals uniquely.
The behavior of the reconstruction in the presence of noise has
been also investigated.

The main type of signals for which explicit reconstruction schemes
have been proposed include linear combinations of
$\delta$-functions and their derivatives, splines, and piecewise
polynomials. In spite of a somewhat different setting of the
problem, the reconstruction schemes turn out to be mathematically
similar to the ones presented below. In fact, moments enter, as an
intermediate step, the reconstruction procedure in \cite{Vet5},
and systems very similar to (2.7) and (2.9) below explicitly
appear in \cite{Vet2,Vet4,Vet5}. It is a remarkable fact (although
traced at least to \cite{Pro}) that exactly the same systems arise
in exponential approximation (\cite{Hil}), in reconstruction of
plane polygons (\cite{Mil1,Mil2,Mil3}, see Section 3.1 below), in
reconstruction of quadrature domains (\cite{Put1}, see Section 3.2
below), in Pad\'e approximations, and in many other problems.

In \cite{Vet3} the approach of \cite{Vet2,Vet4,Vet5} is extended
to some classes of parametric non-bandlimited two-dimensional
signals. This includes linear combinations of $2$D
$\delta$-functions, lines, and polygons. Notice that the first
problem, in its complex setting (where we consider as the allowed
measurements only the {\it complex} moments $\mu_k(f)=\int\int z^k
f(x,y) dxdy$) leads once more to a complex system (2.7).

\subsection{Polynomials}

Let $P(x)$ be a polynomial of degree $d$, $P(x)=\sum_{j=0}^d
a_jx^j$. For the $k$-th moment $m_k(P)$ we have \be
m_k(P)=\int_0^1 \sum_{j=0}^d a_jx^{j+k} dx =\sum_{j=0}^d {a_j\over
{j+k+1}}=\sum_{j=0}^d h_{kj}a_j\ ,\ee if we put $h_{kj}={1\over
{j+k+1}}$. Now let $a$ denote the column-vector of the
coefficients $a_j$ of the polynomial $P(x)$ and let $m$ denote the
column-vector of the moments $m_0(P),\dots,m_d(P)$. We get the
following linear system: \be Ha=m, \ H=(h_{kj}).\ee Notice that
the matrix $H$ is a Hankel matrix: the rows of this matrix are
obtained by the shifts of its first row.  More specifically, the
matrix $H$ belongs to the class of Hilbert-type matrices (see
\cite{Kaly}). In particular, its determinant is nonzero, and
system (2.2) has unique solution. Therefore, we have \bp A
polynomial $P(x)$ of degree $d$ can be uniquely reconstructed from
its first $d+1$ moments $m_0(P),\dots,m_d(P)$, via solving system
(2.3). \ep Notice, however, that the smallest eigenvalue
$\lambda_{min}(H)$ behaves asymptotically for $d \rightarrow
\infty$ as follows:
$$\lambda_{min}(H)=K\sqrt d \rho^{-4(d+1)}(1+o(1)),$$ where
$K=8\pi \sqrt {2\pi}2^{1\over 4}$ and $\rho=1+\sqrt 2$.
(\cite{Kaly}). Therefore, the inversion of the matrix $H$ becomes
problematic for large $d$.

Notice also that for each fixed polynomial $P(x)$ expression (2.1)
defines $m_k(P)$ as a rational function of $k$.

\medskip

As for the moment generating functions, we have \bp
$I_P(z)=-{1\over z}\log(1-{1\over z})+\hat P({1\over z})$, with
$\hat P(s)$ a polynomial of degree $d-1$ in $s$.\ep \pr We have
$P(t)=\tilde P(t)(t-{1\over z}) + P({1\over z})$ where $\tilde
P(t)$ is a polynomial of degree $d-1$ in $t$ with the coefficients
- polynomials of degree $d-1$ in ${1\over z}$. Hence
$$I_P(z)=\int_0^1 {{P(t)dt}\over {1-zt}}=-{1\over z}\int_0^1
{{P({1\over z})dt}\over {t-{1\over z}}}-\int_0^1\tilde P(t)dt.$$
Integrating from $0$ to $1$ now provides the required expression.

\subsection{Rational functions}

Let $R(x)$ be a rational function of degree $d$,
$R(x)={{P(x)}\over {Q(x)}}, \ deg \ Q=d, \ deg \ P \leq d-1$ (we
assume that $R$ does not have a ``polynomial part"). Thus
$$P(x)=\sum_{j=0}^{d-1} a_jx^j, \ Q(x)=\sum_{j=0}^d b_jx^j.$$ We have
$$P(x)=Q(x)R(x)=\sum_{j=0}^d b_jx^jR(x).$$ Hence $$ m_k(P)=
\sum_{j=0}^d b_jm_{k+j}(R), \ k=0,1,\dots,$$ and using our
notations from Section 2.1 above we finally get a system for the
unknowns $a_j, \ b_j$ \be \sum_{j=0}^{d-1} h_{kj}a_j=\sum_{j=0}^d
m_{k+j}(R)b_j, \ k=0,1,\dots,2d,\ee where, as above,
$h_{kj}={1\over {j+k+1}}$. We do not analyze here the solvability
conditions for (2.4) (compare, however, Lemma 2 in \cite{Kis1}).
Let us notice also that the counting of the sign changes as in
Section 2.6 below shows that a rational function $R(x)$ of degree
$d$ can be uniquely reconstructed from its first $4d$ moments
$m_0(R),\dots,m_{4d}(R)$.

\smallskip

To compute the moment generating function $I_R(z)$ let us assume
that the roots $\alpha_1,\dots,\alpha_d$ of $Q$ are all distinct.
Then $R(t)=\sum^d_{i=1}{{A_i}\over {t-\alpha_i}}$ and denoting
${1\over z}$ by $w$ we get ${{R(t)}\over
{t-w}}=\sum^d_{i=1}{{A_i}\over
{(t-\alpha_i)(t-w)}}=\sum^d_{i=1}A_i({1\over{(\alpha_i-w)(t-\alpha_i)}}
-{1\over{(\alpha_i-w)(t-w)}}).$ Transforming integral (2.1) as in
the proof of Proposition 2.2, and integrating we finally get \bp
The moment generating function $I_R(z)$ of a rational function
$R(x)$ is given by \be
I_R(z)=-w\sum^d_{i=1}{A_i\over{\alpha_i-w}}[\log({{1-\alpha_i}\over
{\alpha_i}})-\log({{w-1}\over w})], \ w={1\over z}.\ee \ep

\subsection{Linear combination of $\delta$-functions}

Let $g(x)=\Sigma_{i=1}^n A_i\delta(x-x_i)$. For this function we
have \be m_k(g)=\int_0^1 x^k \Sigma_{i=1}^n A_i\delta(x-x_i) dx =
\Sigma_{i=1}^n A_i x_i^k.\ee So assuming that we know the moments
$m_k(g)=\alpha_k, \ k=0,1,\dots, 2n-1,$ we obtain the following
system of equations for the parameters $A_i$ and $x_i, \
i=1,\dots,n,$ of the function $g$: \be \Sigma_{i=1}^n A_i x_i^k =
\alpha_k, \ k=0,1,\dots, 2n-1.\ee Notice that system (2.7) is
linear with respect to the parameters $A_i$ and non-linear with
respect to the parameters $x_i$.

\smallskip

System (2.7) appears in many mathematical and applied problems.
First of all, if we want to approximate a given function $f(x)$ by
an exponential sum $$ f(x)\approx C_1e^{a_1x}+C_2e^{a_2x}+\dots+
C_ne^{a_nx},$$ then the coefficients $C_i$ and the values
$\mu_i=e^{a_i}$ satisfy a system of the form (2.7) with the
right-hand side (the ``measurements") being the values of $f(x)$
at the integer points $x=1,2,\dots$. (see \cite{Hil}, Section
4.9). The method of solution of (2.7) which we give below, is
usually called Prony's method (\cite{Pro}).

On the other hand, system (2.7), recurrence (2.9) and system
(2.10) below form one of the central objects in Pad\'e
approximation: see, in particular, \cite{Nik.Sor} and references
there.

System (2.7) appears also in error correction codes, in array
processing (estimating the direction of signal arrival) and in
other applications in Signal Processing (see, for example,
\cite{Mil3,Vet5} and references there).

In \cite{Mil1,Mil2,Mil3} system (2.7) appears in reconstruction of
plane polygons from their complex moments. These results are
shortly described in Section 3.2 below.

This system appears also in some perturbation problems in
nonlinear model estimation.

\medskip

We give now a sketch of the proof of solvability of (2.7) and of
the solution method, which is, essentially, the Prony's one. We
follow the lines of \cite{Mil3,Kis2}. See also a literature on
Pad\'e approximation, in particular, \cite{Nik.Sor} and references
there.

\bt A linear combination $g(x)$ of $n$ $\delta$-functions can be
uniquely reconstructed from its first $2n-1$ moments
$m_0(g),\dots,m_{2n-1}(g)$, via solving system (2.7). \et \pr
Representation (2.6) of the moments immediately implies the
following result for the moments generating function $I_g(z)$: \bp
For $g(x)=\Sigma_{i=1}^n A_i\delta(x-x_i)$, the moments generating
function $I_g(z)$ is a rational function with the poles at $x_i$
and with the residues at these poles $A_i$: \be I(z)=
\Sigma_{i=1}^n {A_i\over {1-zx_i}}.\ee \ep We see that the
function $I(z)$ encodes the solution of system (2.7). So to solve
this system it remains to find explicitly the rational function
$I(z)$ from the first $2n$ its Taylor coefficients
$\alpha_0,\dots,\alpha_{2n-1}$. This is, essentially, the problem
of Pad\'e approximation (\cite{Nik.Sor}).

\medskip

Now we use the fact that the Taylor coefficients of a rational
function of degree $n$ satisfy a linear recurrence relation of the
form \be m_{r+n}=\Sigma_{j=0}^{n-1} C_j m_{r+j}, \ r=0,1,\dots.\ee
Since we know the first $2n$ Taylor coefficients
$\alpha_0,\dots,\alpha_{2n-1}$, we can write a {\it linear} system
on the unknown recursion coefficients $C_l$: \be
\Sigma_{j=0}^{n-1} C_j \alpha_{j+r}=\alpha_{n+r}, \
r=0,1,\dots,n-1.\ee Solving linear system (2.10) with respect to
the recurrence coefficients $C_j$ we find them explicitly. For a
solvability of (2.10) see \cite{Hil,Nik.Sor,Mil3,Kis2}. Now the
recurrence relation (2.10) with known coefficients $C_l$ and known
initial moments allows us to easily reconstruct the generating
function $I_g(z)$ and hence to solve (2.7).

\noindent{\bf Remark.} Another proof of Theorem 2.1 can be
obtained in lines of the proof of Theorem 2.2 below. Indeed, a
difference of two linear combinations of $n$ $\delta$-functions
can have at most $2n-1$ ``sign changes". Then we apply Lemma 2.2.

\subsection{Piecewise-solutions of linear ODE's}

In this paper we do not consider separately the case of
piecewise-polynomials. See \cite{Vet2} where a method for
reconstruction of piecewise-polynomials from samplings is
suggested (which starts with a reconstruction of linear
combinations of $\delta$-functions and of their derivatives).
Instead we consider, as a natural generalization of
piecewise-polynomials, the class $A_D$ of piecewise-analytic
functions, each piece satisfying a fixed linear differential
operator $D$ with rational coefficients. Such functions are
usually called ``L-splines" (see \cite{Sch1,Sch2} and references
there). For piecewise-polynomials of degree $d$ we have
$D={{d^{d+1}}\over {dx^{d+1}}}$. Notice that Vetterli's method
(\cite{Vet2}) can be extended also to our class $A_D$. However, in
the present paper we stress another approach to the moment
reconstruction problem for the class $A_D$. It is presented in the
accompanying paper \cite{Kis1} (this volume), while here we
provide a necessary background.

\medskip

Consider the equation \be Dy= y^{(k)}+a_{k-1}(x)y^{(k-1)}+ \dots +
a_1(x)y'+a_0(x)y=0 \ee with the coefficients $a_{k-1}(x), \dots,
a_0(x)$ real-analytic and regular on $[0,1]$. All the solutions of
(2.11) on $[0,1]$ form a linear space $L_D$ with the basis
$y_1(x),\dots, y_k(x)$ being the fundamental set of solutions of
(2.11). For $D={{d^k}\over {dx^k}}$ the space $L_D$ consists of
all the polynomials of degree at most $k-1$, and we can take
$\{y_1(x),\dots, y_k(x)\}=\{1,x,x^2,\dots,x^{k-1}\}$.

\medskip

Now we consider the class $A_D$ of all the piecewise-continuous
functions $g(x)$ on $[0,1]$ with the jumps at $x_1,\dots,x_n\in
[0,1]$, such that on each continuity interval
$\Delta_i=[x_i,x_{i+1}]$ the function $g(x)$ satisfies $Dg=0$. We
extend $g(x)$ by the identical zero outside the interval $[0,1]$.

We can represent $g(x)$ on the intervals $\Delta_i$ in a
``polynomial form": $g(x)=\sum_{j=1}^k \alpha_{ij}y_j(x),$ where
$y_1(x),\dots, y_k(x)$ is the fundamental set of solutions of
(2.11). Alternatively, we can parametrize $g(x)$ on the intervals
$\Delta_i$ by its initial data at the point $x_i$. We can further
define ``splines" of a prescribed smoothness in $A_D$. The
constructions of \cite{Vet2} can be extended to this case.

\medskip

While till this point we could restrict our presentation to the
real domain, in what follows it will be necessary to extend the
consideration to the complex plane.

\medskip

First we recall shortly some classical facts related to the
structure of linear differential equations in the complex domain
(see, for example, \cite{pry,ry} where these fact are presented in
a form convenient for our applications).

\medskip

Consider the equation \be Dy=y^{(k)}+a_{k-1}(x)y^{(k-1)}+ \dots +
a_1(x)y'+a_0(x)y=0 \ee with the coefficients $a_{k-1}(x), \dots,
a_0(x)$ regular and {\it univalued} in the complex domain $\O=
{\mathbb C} \setminus \{x_0,\dots, x_m\}$. We do not specify at
this stage the character of possible singularities of $a_j(x)$ at
the points $x_0,\dots, x_m$.

The following proposition (see, for example, \cite{ry})
characterizes multivalued analytic functions which are solutions
of a certain equation of the form (2.12):

\bp Any solution $y(x)$ of (2.12) is a regular multivalued
function in $\O$, satisfying the following additional property
(F): For any point $w\in \O$ the linear subspace $L_w$ spanned by
all the branches of $y(x)$ at $w$ in the space ${\cal O}(w)$ of
all the analytic germs at $w$, has dimension at most $k$.

Any regular multivalued function $v(x)$ in $\O$ with the property
(F) satisfies a certain equation of the form (2.12) of order at
most $k$ with all the coefficients regular and univalued in the
domain $\O$. \ep Let us remind that for a given function $g(x)$ on
$[0,1]$ the moment generating function $I_g(z)=\sum_{k=0}^{\infty}
m_k(g) z^k$ is given by the Cauchy-type integral \be
I_g(z)=\int_0^1 {{g(t)dt}\over {1-zt}}=w\int_0^1 {{g(t)dt}\over
{w-t}}, \ \ w={1\over z}.\ee Now one of the basic classical facts
about Cauchy-type integrals is that if $g$ (on each its continuity
interval) satisfies a certain equation of the form (2.12) then
$I_g(z)$ satisfies another equation of this form. A proof (in a
specific case which we need in the present paper) can be found in
\cite{pry,ry}. In these papers also specific ramification
properties of $I_g(z)$ are studied  for $g$ algebraic.

\medskip

Now, in the accompanying paper \cite{Kis1} (this volume) the
functions $g(x)$ from the class $A_D$ are considered. A {\it
non-homogeneous} equation of the form (2.12) for $I_g(z)$ is
presented explicitly, and on this base a reconstruction procedure
is suggested.

\subsection{Piecewise-algebraic functions}

Exact reconstruction of piecewise-algebraic ($=$ semi-algebraic)
functions can be considered as one of the ultimate goals of our
approach. If we extend this class $SA$ to
$SA(\psi_1,\dots,\psi_l)$, adding a finite number of fixed
``models" $\psi_1,\dots,\psi_l$ and allowing for all the
elementary operations and for solving equations, we shall probably
cover all the examples of interest. In particular, such extensions
include linear combinations of shifts and dilations of
$\psi_1,\dots,\psi_l$ - an important class appearing in
reconstruction of signals with finite innovation rate
(\cite{Vet2,Vet4,Vet5}), and in wavelets theory. Extensions of
this sort are also closely related to what appears in theory of
$o$-minimal structures - see, for example, \cite{vdd}. Because of
the ``finiteness results" in this theory we can hope that the
``finite moments determinacy" of semi-algebraic functions (Theorem
2.2 below) can be extended to at least some important classes
$SA(\psi_1,\dots,\psi_l)$.

\medskip

Let us remind that $g(x)$ is an algebraic function (as usual,
restricted to $[0,1]$) if $y=g(x)$ satisfies an equation \be
a_n(x)y^n+a_{n-1}(x)y^{n-1}+\dots + a_1(x)y+a_0(x)=0,\ee where
$a_n(x),\dots,a_0(x)$ are polynomials in $x$ of degree $m$.
$d=m+n$ is, by definition, the degree $\deg g$ of $g$.

We shall need the following simple properties of algebraic
functions:

\smallskip

1. The number of zeroes of an algebraic function $g(x)$ defined by
(2.14) does not exceed $m$ (and so it does not exceed its degree
$\deg g=m+n$).

\smallskip

2. A sum $g(x)=g_1(x)+g_2(x)$ of two algebraic functions of
degrees $d_1$ and $d_2$ is an algebraic function, with the degree
$\deg g \leq \eta(d_1,d_2)$.

\medskip

We consider piecewise-algebraic functions on $[0,1]$. Let such a
function $g(x)$ be represented by the algebraic functions $g_q(x)$
of the degrees $d_q$, respectively, on the intervals
$\Delta_q=[x_q,x_{q+1}], \ q=0,...,r,$ of the partition of $[0,1]$
by $x_0=0<x_1<\dots<x_r<x_{r+1}=1$. We define the {\it
combinatorial complexity, (or the degree)} $\sigma(g)$ of $g$ as
follows: \bd (See \cite{Yom1,Yom2}). The combinatorial complexity
$\sigma(g)$ is the sum $\sum^r_{q=1}d_q +r.$ \ed The specific
choice of this expression is motivated by the following simple
observation: {\it the number of sign changes of a
piecewise-algebraic function $g$ on $[0,1]$ does not exceed
$\sigma(g)$.} This follows directly from property (1) above.

We need also the following lemma: \bl Let $g_1,g_2$ be
piecewise-algebraic functions with $\sigma(g_j)\leq d, \ j=1,2$.
Then for $g=g_1\pm g_2$ the combinatorial complexity $\sigma(g)$
satisfies $\sigma(g) \leq  \kappa(d)=2d(\eta(d,d)+1),$ where
$\eta(d,d)$ is given by property (2) above. \el\pr $g$ has at most
$2d$ jumps, and on each continuity interval its degree is bounded
by $\eta(d,d)$.

\smallskip

Now we can show that piecewise-algebraic functions are uniquely
defined by their few moments. We do not touch in this stage the
question of {\it how such a function can be actually reconstructed
from the moments data}, postponing this problem till Section
2.6.1. \bt A piecewise-algebraic function of a combinatorial
complexity $d$ is uniquely defined by its first $\kappa(d)$
moments.\et \pr Assume, in contrary to the statement of the
theorem, that there are functions $g_1$ and $g_2$ of complexity at
most $d$, with exactly the same moments up to order $s=\kappa(d)$.
Hence for the difference $g=g_2-g_1\ne 0$ we have the vanishing of
the moments up to $s$: $m_j(g)=0, \ j=0,1,\dots,s$. By Lemma 2.1
we have for the combinatorial complexity of $g$ the bound
$\sigma(g)\leq s$. Consequently, the number of sign changes of $g$
does not exceed $s$. The next trick comes from the classical
moment theory: \bl If the number of the sign changes and zeroes of
$g(x)\ne 0$ does not exceed $s$ then some of its first $s$ moments
$m_j(g), \ j=0,1,\dots,s$ do not vanish.\el \pr We can assume that
$g$ changes its sign at certain points $t_1,\dots,t_l, \ l \leq
s$, and preserves the sine between these points. Let us construct
a polynomial $Q(t)$ of degree $l$ with exactly the same sign
pattern as $g$: $Q(t)=\pm (x-t_1)(x-t_2)\cdots (x-t_l)$. Write $Q$
as $Q(x)=\sum^l_1 \alpha_j x^j$. We have $g(x)Q(x) >0$ everywhere,
besides $t_1,\dots,t_l$ and possibly some other isolated points.
Therefore $\int_0^1 g(x)Q(x)
>0$. On the other hand, this integral can be expressed as a linear
combination of the moments: $\int_0^1 g(x)Q(x)=\sum^l_1 \alpha_j
\int_0^1 x^j g(x)dx = \sum^l_1 \alpha_j m_j(g).$ Hence some of the
moments of $g$ up to $l\leq s$-th  do not vanish. This proves
Lemma 2.2. To complete the proof of Theorem 2.2 it remains to
notice that the difference $g=g_2-g_1\ne 0$ on at least one of its
continuity intervals.

\subsubsection{Explicit moment inversion for algebraic functions}

As far as an explicit inversion of the moment transform of
algebraic functions is concerned, we are not aware of any general
approach to this problem. Piecewise-algebraic functions belong to
the class $A_D$, as defined in Section 2.3 above. However, the
problem is that we do not know a priori the differential operator
$D$ which annihilates a given algebraic function $g$. (The form of
D is known, but not the coefficients of the rational entries of
$D$). This fact seems to prevent a direct application of the
method of \cite{Kis1} to piecewise-algebraic functions.

\medskip

Let us analyze in more detail one special case. Assume that the
algebraic curve $y=g(x)$ is a rational one. This means that it
allows for a rational parametrization \be x=P(t), \ y=Q(t).\ee The
moments $m_k(g)$ given by $ m_k(g)=\int_0^1 x^k g(x)dx, \
k=0,1,\dots,$ now can be expressed as \be m_k(g)=\int_a^b
P^k(t)Q(t)p(t)dt,\ee where $p$ denotes the derivative of $P$ and
$0=P(a), \ 1=P(b)$. Moments of this form naturally appear in a
relation with some classical problems in Qualitative Theory of
ODE's - see \cite{bry,bfy,bry1,chr}, \cite{pak1}-\cite{pry}.

\medskip

Our problem can be reformulated now as the problem of explicitly
finding $P$ and $Q$ from knowing a certain number of the moments
$m_k$ in (2.16).

\smallskip

Of course, in general we cannot expect this system of nonlinear
equations to have a unique solution. Indeed, while the function
$y=g(x)$ is determined by its moments in a unique way, the {\it
rational parametrization $P,Q$} of this curve in general is not
unique. In particular, let $W(t)$ be a rational function
satisfying $W(0)=0, \ W(1)=1$. Substituting $W(t)$ into $P$ and
$Q$ we get another rational parametrization of our curve: \be
x=\hat P(t), \ y=\hat Q(t), \ \ with \ \ \hat P(t)=P(W(t)), \ \hat
Q(t)=Q(W(t)).\ee Consequently, we can ask the following question:

\medskip

{\it Are all the solutions of (2.16) related one to another via a
composition transform (2.17)?}

\medskip

If the answer to this question is positive, we can restrict our
parametrizations $P,Q$ to be ``mutually prime in composition
sense" (see \cite{rit}) and thus to obtain uniqueness of the
reconstruction.

More generally, the ``inversion problem" for system (2.16) is:

\medskip

{\it To characterize all the solutions of system (2.16) and to
provide an effective way to find these solutions}.

\medskip

\noindent A special case of the inversion problem, in which
definite results have been recently obtained, is the ``Moment
vanishing problem":

\medskip

{\it To characterize all the pairs $P,Q$ for which the moments
$m_k$ defined by (2.16) vanish}.

\medskip

The moment vanishing problem plays a central role in study of the
center conditions for the Abel differential equation (see
\cite{bry,bfy,chr}, \cite{pak1}-\cite{pry}). In fact, it provides
an infinitesimal version of the Poincer\'e Center-Focus problem
for the Abel equation. In spite of a very classical setting (we
ask for conditions of orthogonality of $pQ$ to all the powers of
$P$!) this problem has been solved (for $P$ and $Q$ polynomials)
only very recently (\cite{pak2}). Let us describe the solution.

We say that $P$ and $Q$ satisfy a ``composition condition" if
there are polynomials $\tilde P(w)$ and $\tilde Q(w)$, and a
polynomial $W(x)$, satisfying $W(0)=W(1),$ such that \be
P(x)=\tilde P(W(x)), \ Q(x)=\tilde Q(W(x)).\ee Composition
condition (2.18) can be easily shown to imply the vanishing of all
the moments (2.16). In many cases it is also a necessary one, but
not always. The examples of $P,Q$ annihilating the moments (2.16)
but not satisfying (2.18) can be obtained as follows (see
\cite{pak1}): if $P$ has two right composition factors $W_1(x)$
and $W_2(x)$, then $P$ and $Q=W_1+W_2$ will annihilate the moments
(2.16) because of a linearity with respect to $Q$. For some $P$ we
can find $W_1$ and $W_2$ which are mutually prime in composition
algebra (see \cite{rit}). Then typically $P$ and $Q=W_1+W_2$ will
have no common right composition factors (\cite{pak1}). The result
of \cite{pak3} claims that this is essentially the only
possibility: \bt (\cite{pak3}) All the moments (2.16) vanish if
and only if $Q$ is a sum of $Q_j, \ j=1,\dots,l,$ such that for
each $j$ the polynomials $P$ and $Q_j$ satisfy composition
condition (2.18).\et One can expect that the methods developed in
\cite{bry,bfy,chr,pak3}, \cite{pak1}-\cite{pry} can help in
further analyzing the reconstruction problem for semi-algebraic
functions in one and more variables. See, in particular, Section
3.3. below.

\section{Functions of two variables}
\setcounter{equation}{0}

Also in two dimensions exact reconstruction of semi-algebraic
functions (and of their extension to $SA(\psi_1,\dots,\psi_l)$)
can be considered as one of the ultimate goals of our approach.

\subsection{Reconstruction of polygons from complex moments}

In \cite{Mil3,Mil2,Mil1} the problem of reconstruction of a planar
polygon from its complex moments is considered. The complex
moments of a function $f(x,y)$ are defined as \be
\mu_k(f)=\int\int z^k f(x,y) dx dy, \ k=0,1,\dots, \ z=x+iy. \ee
Complex moments can be expressed as certain specific linear
combinations of the real double moments $m_{kl}(f)$.

For a plane subset $A$ its complex moments $\mu_k(A)$ are defined
by $\mu_k(A)=\mu_k(\chi_A),$ where $\chi_A$ is the characteristic
function of $A$.

\smallskip

Let $P$ be a closed $n$-sided planar polygon with the vertices
$z_i, \ i=1,\dots,n$. The reconstruction method of \cite{Mil3} is
based on the following result of \cite{Dav}: \bt There exist a set
of $n$ coefficients $a_i, \ i=1,\dots,n$, depending only on the
vertices $z_i$, such that for any analytic function $\phi(z)$ on
$P$ we have $$\int\int_P \phi''(z) dxdy=\sum_{i=1}^n
a_i\phi(z_i).$$ The coefficients $a_j, \ j=1,\dots,n$ are given as
$a_j={1\over 2}({{\bar z_{j-1}-\bar z_j}\over {z_{j-1}-z_j}} -
{{\bar z_{j}-\bar z_{j+1}}\over {z_{j}-z_{j+1}}})$.\et Applying
this formula to $\phi(z)=z^k$ we get \be
k(k-1)\mu_{k-2}(\chi_P)=\sum_{i=1}^n a_i z_i^k, \ k=0,1,\dots,\ee
where we put $\mu_{-2}=\mu_{-1}=0$. So on the left-hand side we
have shifted moments of $P$.

\smallskip

If we ignore the fact that $a_j$ can be expressed through $z_i$
and consider both $a_j$ and $z_i$ as unknowns, we get from (3.2) a
system of equations \be \sum_{i=1}^n a_i z_i^k=\nu_k, \
k=0,1,\dots,\ee where $\nu_k$ denotes the ``measurement"
$k(k-1)\mu_{k-2}(P)$. System (3.3) is identical to system (2.7)
which appears in reconstruction of linear combination of
$\delta$-functions. One of the solution methods suggested in
\cite{Mil3} is the Prony method described in Section 2.4 above.
Another approach is based on matrix pencils. In \cite{Mil2,Mil1}
an important question is investigated of polygon reconstruction
from noisy data.

\subsection{Quadrature domains}

We introduce, following \cite{Put1}, a slightly different sequence
of double moments: for a function $g(z)=g(x+iy)$ the moments
$\tilde m_{kl}(g)$ are defined by \be \tilde m_{kl}(g)=\int\int
z^k \bar z^l g(z) dx dy, \ k,l \in {\mathbb N}.\ee One defines the
moment generating function $I_g(v,w)=\sum_{k,l=0}^{\infty}\tilde
m_{kl}(g)v^kw^l$ and the ``exponential transform" $$ \tilde
I_g(v,w)= 1-\exp(-{1\over \pi}I_g(v,w))=$$ $$ =\exp(-{1\over
\pi}\int\int_\Omega{{g(z)dxdy}\over {(z-v)(\bar
z-w)}}):=\sum_{k,l=0}^{\infty}b_{kl}(g)v^kw^l.$$ Now (classical)
quadrature domains in ${\mathbb C}$ are defined as follows: \bd A
quadrature domain $\Omega \subset {\mathbb C}$ is a bounded domain
with the property that there exist points $z_1,\dots,z_m \in
\Omega$ and coefficients $c_{ij}, \ i=1,\dots,m, \ j=0,\dots,\
s_i-1$, so that for all analytic integrable functions $f(z)$ in
$\Omega$ we have \be \int\int_\Omega f(x+iy)
dxdy=\sum_{i=1}^m\sum_{j=0}^{s_i-1}c_{ij}f^{(j)}(z_i).\ee
$N=s_1+\dots+s_m$ is called the order of the quadrature domain
$\Omega$. \ed The simplest example is provided by the disk
$D_R(0)$ of radius $R$ centered at $0\in {\mathbb C}$:
$\int\int_{D_R(0)} f(x+iy) dxdy=\pi R^2f(0)$. The results of Davis
(\cite{Dav}; Theorem 3.1 above) give another example in this
spirit.

\smallskip

The following result (\cite{Put2},\cite{Put1}, Theorem 3.1)
provides a necessary and sufficient condition for $\Omega \subset
{\mathbb C}$ to be a quadrature domain: let $\tilde
I_\Omega(v,w)=\tilde I_{\chi_\Omega}(v,w)$ be the exponential
transform of $\Omega$. \bt $\Omega$ is a quadrature domain if and
only if there exists a polynomial $p(z)$ with the property that
the function $\tilde q(z,\bar w)=p(z)\bar p(w)\tilde
I_\Omega(z,\bar w)$ is a polynomial at infinity (denoted by
$q(z,\bar w)$). In that case, by choosing $p(z)$ of minimal
degree, the domain $\Omega$ is given by $\Omega=\{z\in {\mathbb
C}, \ q(z,\bar z)<0\}.$ Moreover, the polynomial $p(z)$ in this
case is given by $p(z)=\prod_{i=1}^m(z-z_i)^{s_i}$, where $z_i$
are the quadrature nodes of $\Omega$.\et Now, the algorithm in
\cite{Put1} for reconstruction of a quadrature domain from its
moments consists of the following steps:

\medskip

1. Given the moments $\tilde m_{kl}(\Omega)=\tilde
m_{kl}(\chi_\Omega)$ up to a certain order, construct the
(truncated) exponential transform $\tilde I(v,w) =
\sum_{k,l=0}^{\infty}b_{kl}v^kw^l$.

\smallskip

2. Identify the minimal integer $N$ such that
$\det(b_{k,l})_{k,l=0}^N=0$. Then there are coefficients
$\alpha_j, \ j=0,\dots,N-1,$ such that for $B=(b_{k,l})_{k,l=0}^N$
and $\alpha =(\alpha_1,\dots,\alpha_{N-1},1)^T$ we have \be
B\alpha=0.\ee We solve this system with respect to $\alpha$. Then
the polynomial $p(z)$ defined above is given by
$p(z)=z^N+\alpha_{N-1}z^{N-1}+\dots+\alpha_0.$

\smallskip

3. Construct the function $$R_\Omega(z,\bar w)=p(z)\bar p(w)
exp(-{1\over \pi}\sum_{k,l=0}^{N-1}\tilde m_{kl}(\Omega){1\over
{z^{k+1}}}{1\over {\bar w^{l+1}}})$$ and identify $q(z,\bar w)$ as
the part of $R_\Omega(z,\bar w)$ which does not contain negative
powers of $z$ and $\bar w$. Then the domain $\Omega$ is given by
$\Omega=\{z\in {\mathbb C}, \ q(z,\bar z)<0\}.$

\medskip

\noindent{\bf Remark.} Let us substitute into the definition of
the quadrature domain (formula (3.5) above) $f(z)=z^k$. Assuming
that all the quadrature nodes $z_i$ are simple, we get for the
complex moments $\tilde m_{k,0}(\Omega)=m_k(\Omega)$ the
expression $$m_k(\Omega)= \sum_{i=1}^m c_i z_i^k,$$ which is
identical to (3.3) in reconstruction of planar polygons. So we can
reconstruct the quadrature nodes $z_i$ and the coefficients $c_i$
from the complex moments only, and we get once more a complex
system which is identical to (2.7). Allowing quadrature nodes
$z_i$ of an arbitrary order, we get a system corresponding to a
linear combination of $\delta$-functions and their derivatives
(compare \cite{Vet2,Vet3}).

Notice that system (3.6) that appears in step 2 of the
reconstruction algorithm above is very similar to system (2.10) in
the solution process of (2.7).

\subsection{$\delta$-functions along algebraic curves}

As we've mentioned above, a natural class of functions $f(x,y)$ of
two variables, for which we can hope for an explicit
reconstruction from a finite number of the moments
$m_{kl}(f)=\int\int x^k y^l f(x,y) dx dy, \ k,l=0,1,\dots,$
consists of semi-algebraic functions. Those are
piecewise-algebraic functions with the continuity pieces bounded
by piecewise-algebraic curves. Among semi-algebraic functions are
piecewise-polynomials with the continuity pieces bounded by spline
curves - a very natural and convenient object in constructive
approximation.

Most of the methods presented in Section 2 for functions of one
variable are applicable also in the case of two variables. In
particular, generalizing the approach of \cite{Vet2} we can
differentiate piecewise-polynomials a sufficient number of times
and finally get a combination of weighted $\delta$ functions along
the partition curves. See also \cite{Vet3}.

In this paper we restrict ourself to a discussion of only one
example. Assume that $f(x,y)$ is a $\delta$-function $\delta_S$
along a rational curve $S$, i.e. for any $\psi(x,y)$ we have
$\int\int f\psi dx dy = \int_{S}\psi(x,y)dx$. Let \be x=P(t), \
y=Q(t), \ t \in [0,1] \ee be a rational parametrization of $S$.
The moments now can be expressed as \be m_{kl}(f)=\int_0^1
P^k(t)Q^l(t)p(t)dt,\ee where $p$ denotes the derivative $P'$ of
$P$. This system is an extension of system (2.16): here we are
allowed to use all the double moments, while in (2.16) only the
moments $m_{k1}$ are available.

\smallskip

Also here we notice that a rational parametrization $P,Q$ of the
curve $S$ in general is not unique: for any rational function
$W(t)$ satisfying $W(0)=0, \ W(1)=1$ we get another rational
parametrization of our curve: \be x=\hat P(t), \ y=\hat Q(t), \ \
with \ \ \hat P(t)=P(W(t)), \ \hat Q(t)=Q(W(t)).\ee Consequently,
we can reiterate the question in Section 2.6 with better chances
for a positive answer:

\medskip

{\it Are all the solutions of (3.9) related one to another via a
composition transform (3.10)?}

\medskip

For the ``Moment vanishing problem" for (3.9) a definite answer
has been obtained in \cite{pry}: {\it composition condition (2.18)
is necessary and sufficient for the moments vanishing.}

\medskip

Let us now assume that the curve $S$ is closed and that it can be
parametrized by $x=P(t), \ y=Q(t)$, with $t$ in the unit sircle
$S^1$. The study of the double moments of this form brings us
naturally to the recent work of G. Henkin \cite{Hen1,Hen2,Hen3}.
Indeed, the vanishing condition for the moments (4.12) is given by
Wermer's theorem (\cite{Wer}): $m_{kl}(f)\equiv 0$ if and only if
$S$ bounds a complex $2$-chain in ${\mathbb C}^2$. See \cite{bry}
for a simple interpretation of Wermer's condition in the case of
rational $P,Q$. In general, if the moments $m_{kl}(f)$ do not
vanish identically, then the local germ of complex analytic curve
$\hat S$ generated by $S$ in ${\mathbb C}^2$ does not ``close up"
inside ${\mathbb C}^2$. G. Henkin's work (\cite{Hen1,Hen2,Hen3}),
in particular, analyzes various possibilities of this sort in
terms of the ``moments generating function". We expect that a
proper interpretation of the results of \cite{Hen1,Hen2,Hen3} can
help also in understanding of the moment inversion problem.

\vskip1cm

\bibliographystyle{amsplain}

\end{document}